\documentclass[12pt]{article}
\usepackage{amsfonts}
\usepackage{amssymb}
\usepackage{amsmath}
\usepackage{epsfig}

\title{The Optimal Twisted Paper Cylinder}
\author{Noah Montgomery and Richard Evan Schwartz \thanks{R.E.S. is supported by
    N.S.F. Grant DMS-2102802 and a Simons Sabbatical Fellowship}}

\newtheorem{theorem}{Theorem}[section]

\newtheorem{lemma}[theorem]{Lemma}

\def\startproof{{\bf {\medskip}{\noindent}Proof: }}

\def\endproof{$\spadesuit$  \newline}

\def\R{\mbox{\boldmath{$R$}}}%
\def\Z{\mbox{\boldmath{$Z$}}}%

\begin{document}
\maketitle

\begin{abstract}
  An embedded twisted paper cylinder of
  aspect ratio $\lambda$ is a smooth isometric
  embedding of a flat $\lambda \times 1$
  cylinder into $\R^3$ such that the
  images of the boundary components
  are linked.  We prove that for such
  an object to exist we must have
  $\lambda>2$ and that this bound is sharp.
  We also show that any sequence of examples
  having aspect ratio converging to $2$
  must converge to a (non-smooth)
  $4$-fold wrapping of a right-angled
  isosceles triangle.
    \end{abstract}

\section{Introduction}

\subsection{Context}

The purpose of this paper is to prove the analogues
for cylinders of the results about Moebius bands
proved in [{\bf S1\/}].

The paper [{\bf S1\/}] deals with the
cutoff value for the aspect ratio of a
rectangle which one can twist in space
to make an embedded paper Moebius band.
The main result of [{\bf S1\/}],
conjectured in 1977 by Ben Halpern and
Charles Weaver [{\bf HW\/}], is 
that a smooth embedded paper Moebius
band must have aspect ratio
greater than $\sqrt 3$.   A secondary
result in [{\bf S1\/}] is that any sequence
of embedded paper Moebius bands with
aspect ratio converging to $\sqrt 3$ must
converge (up to isometries) to a
limiting example called the triangular
Moebius band.   The triangular Moebius band
is a certain $3$-fold wrapping of an
equilateral triangle.

Since [{\bf S1\/}] appeared there has been more
work on related topics:
\begin{itemize}
\item The paper [{\bf S2\/}]  gives an
  explicit estimate for the convergence in
    [{\bf S1\/}]: If $\epsilon<1/324$ and
  a paper Moebius band has aspect ratio
  less than $\sqrt 3 + \epsilon$, then it is
  within $18 \sqrt \epsilon$ of an equilateral
  triangle in the Hausdorff metric.
\item  In [{\bf BrS\/}], Brienne Brown and the second author show that
  one can make a  $3$-twist paper Moebius band having
  aspect ratio $\lambda$ for any $\lambda>3$.   They
  also present  two polygonal
  paper Moebius bands of aspect
  ratio $3$ and conjecture that any
  aspect-ratio minimizing sequence must converge,
  on a subsequence, to one of them.
\item In [{\bf H\/}], Aidan Hennessey showed that
  one can make an embedded
  paper Moebius band or a cylinder with as
  many twists as desired using a rectangle of
  aspect ratio $6$.  In an unpublished observation,
  Jan Neinhaus optimized Hennessey's contruction
  and showed that, for Moebius bands,
  the aspect ratio $3 \sqrt 3 + \epsilon$
  suffices for any $\epsilon>0$.
\end{itemize}

Paper Moebius bands and cylinders
are closely related to {\it folded
  ribbon knots\/}.
Informally, a folded ribbon knot is a planar polygonal
version of a paper Moebius band or cylinder, with
extra combinatorial information keeping track of
``infinitesimal over-crossings and under-crossings''.
See [{\bf DL\/}] and [{\bf D\/}] for a wealth
of information about these.  With a bit of soft
work concerning smooth approximations of folded
ribbon knots (along the lines of \S \ref{exist} below)
our results below resolve the
case $n=1$ of [{\bf DL\/}, Conjecture 39].

\subsection{The Aspect Ratio Bound}

An {\it embedded paper cylinder of aspect ratio $\lambda$\/} is a
$C^{\infty}$, injective, arc-length preserving map
\begin{equation}
  I : \Gamma_{\lambda} \to \R^3,
\end{equation}
where $\Gamma_{\lambda}$ is the flat cylinder of aspect ratio $\lambda$:
\begin{equation}
  \Gamma_{\lambda}=[0,\lambda] \times [0,1]/\sim, \hskip 30 pt
  (0,y) \sim (\lambda,y).
\end{equation}

We call an embedded paper cylinder {\it twisted\/} if the two boundary components
form a nonsplit link.  (This means that the two components cannot be
separated
by a topological $2$-sphere.)
The reason for the name is that if you take a rectangular
strip
of paper, give it a $360$ degree twist, and join the ends together in
space,
then you get a paper cylinder whose boundary components form a
Hopf link, the simplest example of an nonsplit link.

\begin{theorem}
  \label{one}
  An embedded twisted paper cylinder has aspect ratio greater than $2$.
\end{theorem}

We give two proofs of
Theorem \ref{one}.
One is similar to the
proof in [{\bf S1\/}], and the other 
uses some ideas about convex
hulls.  With small changes, our proofs also
work for immersed paper cylinders with linked boundary components.
This is in contrast to the results in [{\bf S1\/}], which
depend crucially on the embedding property.
The main change is that
Lemma \ref{BEND} should be formulated in the immersed case.
See the remark at the end of \S \ref{app}.

We can also weaken the linking hypothesis.
Let $F$ and $G$ be the boundary components of
a paper cylinder $\Omega$.  Our first proof works under
the hypothesis that, for every linear
projection $L$, the two curves
$L(F)$ and $L(G)$ intersect.  Our second proof works
under the hypothesis that $F$ intersects the convex
hull of $G$.  The linking hypothesis implies each of
these other hypotheses.

\subsection{Limiting Behavior}

As we prove below,
an embedded twisted paper cylinder of aspect ratio nearly $2$
must approximate a $4$-fold wrapping
of a right isosceles triangle.  There are essentially
four distinct ways to produce such a wrapping.  Two are shown in
Figure 1.1, and
the other two are their mirror images.  In all four cases  the final result,
shown at right, approximates a right isosceles triangle, and the
boundary components form a Hopf link.

\begin{center}
\resizebox{!}{1.6 in}{\includegraphics{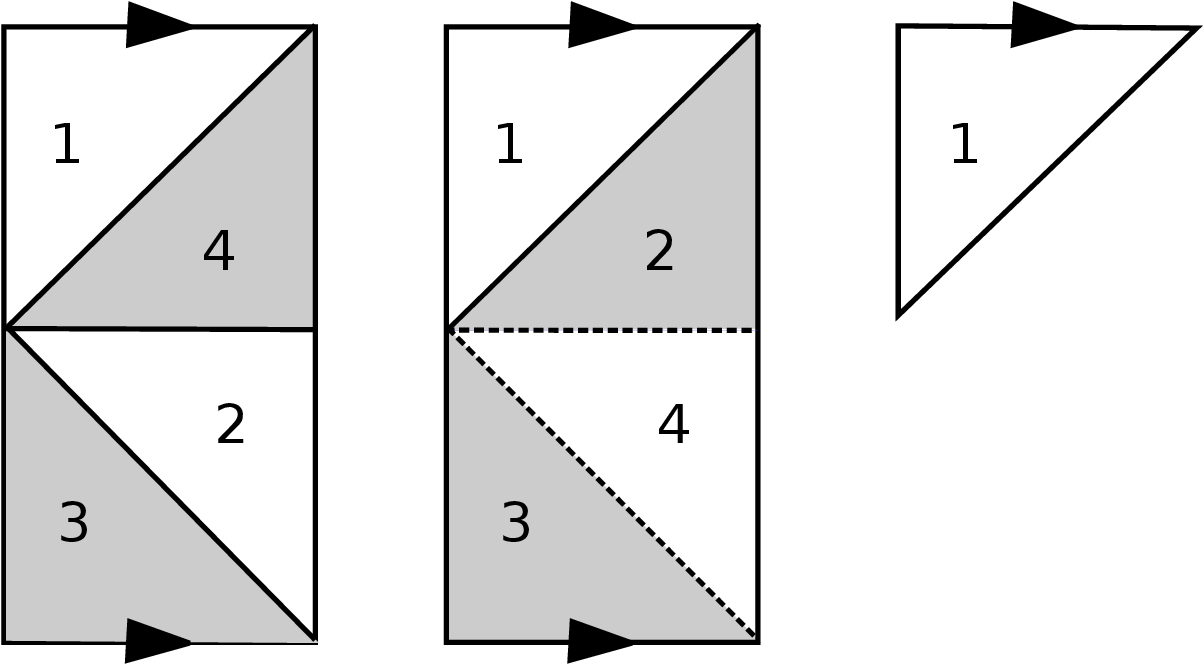}}
\end{center}
\begin{center}
  {\bf Figure 1.1:\/} Two folding patterns and the result of folding them
  \end{center}

  The front sides of the rectangles in Figure 1.1 are divided into $4$
right isosceles
  triangles and marked $1,2,3,4$.  The back sides
 (which are not shown) are unmarked.  We assume
that the
aspect ratio is slightly greater than $2$,
so that the ``segments'' between the triangles have a
tiny bit of thickness,  the folds are more
like tight smooth $U$-turns, and the triangles do not end up in the
same plane after folding.  The edges marked with arrows should be
taped together.  The paper rectangle should be folded down along
the solid segments (so that the marked side of the paper ends up
outside the crease) and up along the dotted segments.   In both case,
the faces end up stacked on top of each other.  If you were
to
poke the stack with a needle you would puncture the faces
in the order $1,2,3,4$.  Within the stack, the white triangles face up and
the
grey triangles face down.

There are other ways to fold up the rectangle to obtain a right isosceles
triangle, but, except for the folding patterns above and their
mirror images, the construction would fail  to produce an embedded
twisted paper cylinder.  Either it would be impossible to tape the edges
without
self-intersection, in violation of the embedding condition, or the
boundary components would fail to be linked.

If the aspect ratio is exactly $2$, all four folding patterns
degenerate into the same continuous and piecewise isometric map
from the cylinder to a right isosceles triangle.   We call this
map the {\it right-isosceles cylinder map\/}.  The image
of the right-isosceles cylinder map, as a set, is just a right
isosceles triangle in space.  However, if we ``remember'' the
folding pattern we could describe the final image as a folded
ribbon knot.
We show that the right-isosceles cylinder
map is the only possible limit of a sequence
of aspect-ratio minimizing embedded twisted paper cylinders.

\begin{theorem}
  \label{three}
  Suppose $I_n: \Gamma_{\lambda_n} \to \R^3$ is a sequence of embedded
  twisted paper cylinders such that $\lambda_n \to 2$.
  Then  $I_n$ converges uniformly,
  up to isometries, to the right-isosceles cylinder map.
\end{theorem}

Corresponding to our two proofs of
Theorem \ref{one}, we give two proofs of
Theorem \ref{three}.   Our two proofs of
Theorem \ref{three} are fairly similar to each other
and indeed have the same endgame.

To show that Theorem \ref{three} is not an empty statement,
we prove that the construction illustrated in Figure 1.1
really does produce (smooth) embedded twisted paper cylinders.
\begin{theorem}
  \label{two}
  For any $\epsilon>0$ there exists an embedded twisted  paper cylinder
  of aspect ratio $2+\epsilon$ whose boundary
  components form a Hopf link.
\end{theorem}
The proof of Theorem \ref{two} is similar in spirit to what
is done in [{\bf HW\/}].  This is just a rigorous working-out
of the description that we gave of the folding
process illustrated by Figure 1.1.

\subsection{Organization}

  This paper is organized as follows.
  In \S 2,  we  give the first proof of
  Theorems \ref{one} and \ref{three}.
  In \S 3,  we  give the second proof.
  In \S 4,  we  prove Theorem \ref{two}.
  In \S 5, an appendix, we prove Lemma \ref{BEND},
  a key ingredient in the proofs.

  \subsection{Acknowledgements}
  
R.E.S.  would like to thank Brienne Brown, Elizabeth Denne,
Eliot Fried, Jeremy Kahn, Curtis McMullen, and
Sergei Tabachnikov for conversations related to this
paper.  R.E.S. also thanks the National
Science Foundation and the Simons Foundation for
their support.

    \section{The First Proof}

\subsection{The Bend Partition}
\label{BEND0}

Let $\Omega=I(\Gamma_{\lambda})$ be an embedded
twisted paper cylinder.
A {\it bend\/} on $\Omega$ is a straight
line segment having one endpoint in each
boundary component of $\Omega$.
In an appendix we
establish the following classic result.
\begin{lemma}
  \label{BEND}
  $\Omega$ has a continuous partition into bends.
\end{lemma}
One property we will frequently use is that bends
have length at least $1$.

The preimage of a bend on $\Omega$, which we
call a {\it prebend\/}, is a line segment
that has one endpoint in each component of
$\partial \Gamma_{\lambda}$.

\subsection{Proof of Theorem \ref{one}}
\label{balance}

Throughout the proof,
$\ell(\cdot)$ denotes arc length.
A pair  $(u,v)$ of distinct bends
on $\Omega$ partitions
$\partial \Omega$ into $4$ arcs, $2$ per component.
We call $(u,v)$ a {\it balanced pair\/}
if these $4$ arcs all have the same length.
The following result is an analogue of
Lemma T in [{\bf S1\/}], but easier to prove.

\begin{lemma}
  $\Omega$ has a balanced pair of bends.
\end{lemma}

\startproof
We identify
the midline $([0,\lambda] \times \{1/2\})/\!\sim$
of $\Gamma_{\lambda}$
with $\R/\lambda \Z$.   We parametrize the
bends by where the corresponding prebends intersect
this midline.   Thus
$\beta_t$ denotes the bend such that
the prebend $\beta'_t$ intersects $\R/\lambda \Z$ at
$t \in \R/\lambda \Z$.
Let $s(t)$ be the
slope of $\beta'_t$.  Let
$f(t)=s(t)-s(t+\lambda/2)$.
By construction
$f(t+\lambda/2)=-f(t)$. Hence, by the
intermediate value theorem, there is some
$t$ such that $f(t)=0$.  But then $s(t)=s(t+\lambda/2)$ and
so $\beta_t$ and $\beta_{t+\lambda/2}$ do the job.
\endproof

\begin{lemma}
  \label{LINE}
  Let $A, B \subset \R^2$ be two segments with
  $\ell(A), \ell(B) \geq 1$.  Let
$A_1,A_2$ be the endpoints of $A$ and
let $B_1,B_2$ be the endpoints of $B$.
Let $C_j$ be an
arc that connects $A_j$ to $B_j$.
If $C_1 \cap C_2 \not =\emptyset$ then
$\ell(C_1)+\ell(C_2) \geq 2$.
\end{lemma}

\startproof
Figure 2.1 shows several possibilities for $A$ and $B$;
the proof is the same in all cases. 
\begin{center}
\resizebox{!}{1.3in}{\includegraphics{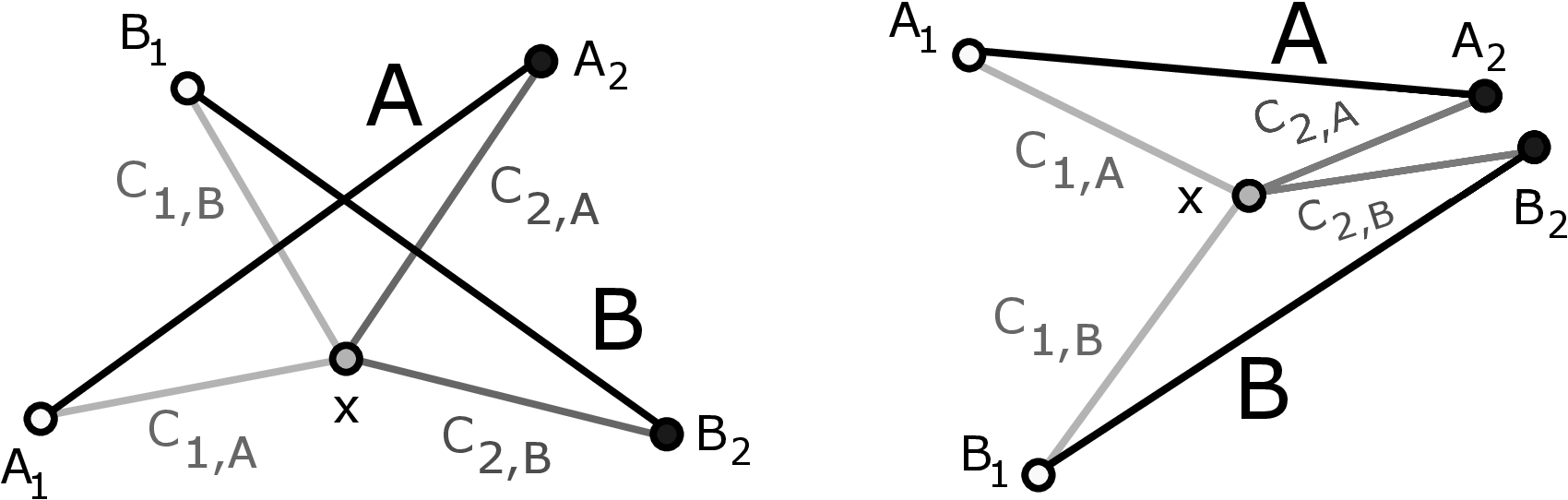}}
\end{center}
\begin{center}
{\bf Figure 2.1:\/} The relevant sets in Lemma \ref{LINE}
\end{center}
A point $x \in C_1 \cap C_2$
divides $C_j$ into paths
$C_{j,A}$ and $C_{j,B}$.  We have
\begin{equation}
  \label{bigon}
\ell(C_1)+\ell(C_2)=
\ell(C_{1,A})+\ell(C_{2,A}) + 
\ell(C_{1,B})+\ell(C_{2,B})
\geq \ell(A)+\ell(B) \geq 2.
\end{equation}
This completes the proof.
\endproof

Let $(u,v)$ be a balanced pair of bends.
We identify $\R^2$ with a plane in $\R^3$.
We move by an isometry so that $u \subset \R^2$ 
and $v$ lies in a plane parallel to $\R^2$.
Let $\Pi$ be the projection from $\R^3$ into $\R^2$.
    Let $A=\Pi(u)=u$ and $B=\Pi(v)$.
    Note that $\ell(A)=\ell(u) \geq 1$ and
      $\ell(B)=\ell(v) \geq 1.$
      Let $F$ and $G$ be the two boundary components of $\Omega$.
      Since these curves are linked,
     the sets $\Pi(F)$ and $\Pi(G)$ intersect.
    In particular there are
    arcs $F^* \subset F$ and $G^* \subset G$, each
    joining an endpoint of $u$ to an endpoint of $v$,
    such that $C_1=\Pi(F^*)$ and $C_2=\Pi(G^*)$ intersect.
    The balance condition and Lemma \ref{LINE} together give:
        \begin{equation}
      \label{threeX}
      \lambda = \ell(F^*)+\ell(G^*) \geq \ell(C_1)+\ell(C_2) \geq 2.
    \end{equation}
    The set of achievable aspect ratios for smooth embedded  paper
    cylinders is open, because we can take any
    example and fatten it a bit.  Hence $\lambda>2$.
    
\vspace{-8 pt}
        \subsection{Proof of Theorem \ref{three}}
        \label{proof31}

        Rather than work with a sequence
         $I_n: \Gamma_{\lambda_n} \to \R^3$ with
         $\lambda_n \to 2$, we will work with
         a single embedded twisted paper cylinder of aspect ratio
         $\lambda=2 + \epsilon$ and examine what happens
         as $\epsilon \to 0$. To save words, we use
        ``nearly''  to mean ``up to an error that tends
          to zero as $\epsilon \to 0$'' and the
          symbol $\approx$ to mean ``nearly equal''.
          Whenever we say two embeddings of curves or surfaces
          are nearly equal, we mean
          with respect to the uniform metric on maps.  Thus,
          our goal
          is to show that $I: \Gamma_{\lambda} \to
          \R^3$ is nearly the right-isosceles cylinder map.
          
          We retain the notation from the preceding subsection.

          Since $\ell(F^*) = \ell(G^*) = \lambda/2 \approx 1$ and
          $\ell(F^*) \geq \ell(C_1)$ and
          $\ell(G^*) \geq \ell(C_2)$,  we see from \eqref{threeX}
         that $\ell(F^*) \approx \ell(C_1)$ and
         $\ell(G^*) \approx \ell(C_2)$.   Since $F^*=\Pi(C_1)$ and
         since $A_1$ is an endpoint of both $F^*$ and $C_1$, we have
          $F^* \approx C_1$.  Likewise $G^* \approx C_2$.  Also
          $u=A=\Pi(u)$ and $v \approx
          B=\Pi(v)$.

          Since $\lambda \approx 2$, all the inequalities in
          \eqref{bigon} are nearly equalities.  In particular
          $\ell(A) \approx \ell(C_{1,A})+\ell(C_{2,A})$.   Hence
          $x$ nearly lies in $A$.  Likewise $x$ nearly lies in $B$.
          Since $\ell(F^*) +  \ell(G^*) \approx 2 \approx \ell(A)+\ell(B)$, we see
          that
          $F^*$ and $G^*$ nearly agree with the
          bigonal paths $A_1,x,B_1$ and $A_2,x,B_2$, respectively.

          Let $\beta$ be a bend whose endpoint in $F^*$ is nearly $x$.
          The endpoint $y$ of $\beta$ in $G$ has distance nearly at
          least $1$ from $x$.
          Since $G$ is a loop of length nearly $2$ which
          contains $y$ and nearly contains $x$, we see that
          $G$ is nearly a unit length bigon having $x$ as an endpoint.
          Likewise, $F$ is nearly a unit length bigon having $x$
          as an endpoint.
          Comparing these descriptions of $F$ and $G$
          with those of their respective subarcs $F^*$ and $G^*$ we
          see that either
          $x \approx A_1 \approx B_2$ or
          $x \approx A_2 \approx B_1$ or $A$ and $B$ nearly overlap on
          a segment whose length is not nearly zero.  The last case is
          impossible because we have just
          shown
          that {\it any\/} point nearly in $A$ and $B$ is nearly a
          common
          endpoint of $F$ and $G$.  So we can assume without
            loss of generality that $x \approx A_2 \approx B_1$.
          Then  $F$ and $G$ respectively are nearly the bigons
              connecting
              the endpoints of $u$ and $v$.
              \newline
              \newline
\noindent
{\bf The Endgame:\/}
We specially highlight this part of the proof because
it is also the endgame for
our second proof of Theorem \ref{three}.

\begin{center}
\resizebox{!}{1.3in}{\includegraphics{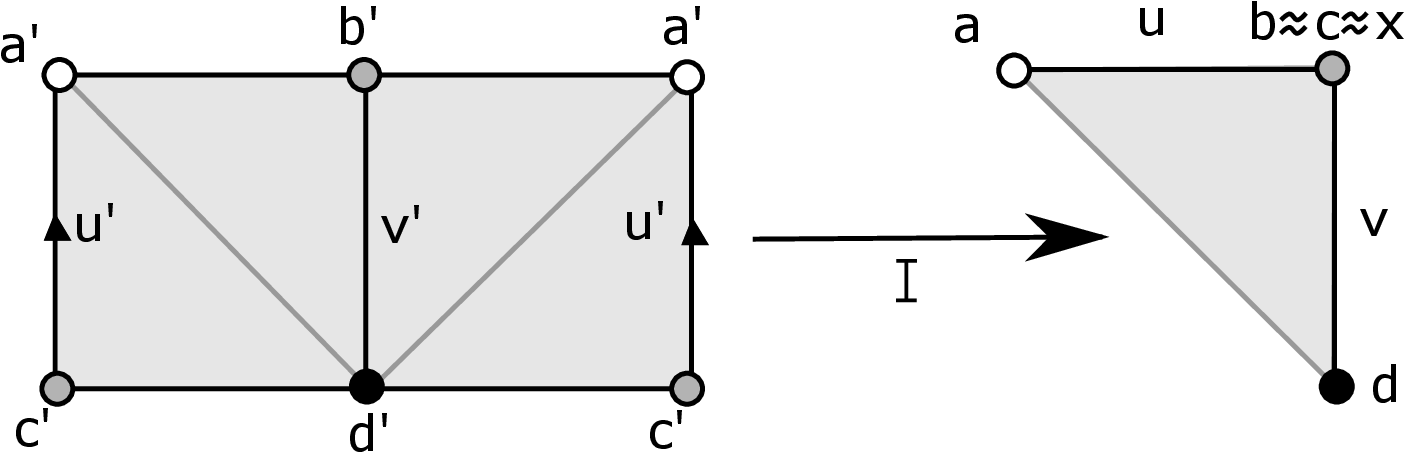}}
\end{center}
\begin{center}
  {\bf Figure 2.2:\/} The prebends and bends and their endpoints
\end{center}

Since $u$ and $v$ have length nearly $1$, the corresponding
prebends $u'$ and $v'$ are nearly perpendicular to
$\partial \Gamma_{\lambda}$ as shown in Figure 2.2.
We label the ends of the prebends and bends as in
Figure 2.2.    Since $\ell(F) \approx 2$ and each of the two arcs of
$F$ connecting $a$ to $b$ has length nearly $1$,
each of these two arcs
has length nearly $1$.  Thus $u'$ and $v'$ nearly divide $\Gamma_{\lambda}$ into
 two unit squares.
\newline
\newline
{\bf Remark:\/} We have refrained from using the fact that
$(u,v)$ is a balanced pair, which would shorten the
argument in the previous paragraph, because we will not have
this situation when we give our second proof of
Theorem \ref{three}.
\newline

\begin{lemma}
  The angle $\theta$ between the vectors $a-c$ and $d-b$ is nearly
  $\pi/2$.
\end{lemma}

\startproof
We use the symbol $s \lesssim t$ to mean
that $\max(0,s-t) \approx 0$.

Since $I$ preserves arc lengths and
$\|a'-d'\| \approx \sqrt 2$ we have
$\|a-d\| \lesssim \sqrt 2$.  This combines with
$b \approx c$ to show that
$\theta \lesssim \pi/2$.

Let $V$ be the ball with $v$ as a diameter.
For any positive $r \leq 1$ we can find a point
$p \in F$ such that $\|c-p\| \approx r$.
Now $p$ is the endpoint of a bend whose
other endpoint $q$ lies in $G$.   But $V$
has diameter nearly $1$ and contains $q$.
Hence $p$ nearly lies outside $V$.
This is equivalent to $\cos(\theta) \lesssim r$.
Since this is true for all $r \in (0,1)$ we have
$\pi/2 \lesssim \theta$.
\endproof

Consider the four triangles into which
$\Gamma_{\lambda}$ is divided in Figure 2.2.
By Lemma 2.4 and the fact that
$\ell(u) $ and $\ell(v)$ are nearly $1$, the embedding
$I$ is nearly a linear isometry on the boundary
of each of these triangles.  Hence $I$ must nearly
be a linear isometry when restricted to each solid
triangle.    We now can conclude that $I$ is nearly
the right-isosceles cylinder map.

    \section{The Second Proof}

\subsection{Proof of Theorem \ref{one}}

Our second proof of Theorem \ref{one} also
relies on Lemma \ref{BEND}
and the terminology established in
\S \ref{BEND0}.

For any subset $S \subset \R^n$ we let
${\rm Hull\/}(S)$ denote the convex hull of $S$,
i.e. the intersection of all convex sets containing $S$.
We will only be working with $n=3$, but the following
well-known lemma works in general.

\begin{lemma}
  \label{hull}
  A subset  $S \subset \R^n$ has the same diameter as its convex hull.
\end{lemma}

\startproof
Let $H={\rm Hull\/}(S)$. Since $S \subset H$, we
have ${\rm diam\/}(S) \leq {\rm diam\/}(H)$.
To show  ${\rm diam\/}(S) \geq {\rm diam\/}(H)$ it
suffices to prove that
for all $p \not = q \in H$ there exist $p',q' \in S$ with
$\|p'-q'\| \geq \|p-q\|$.  Let
$P$ and $Q$ be the disjoint closed halfspaces
bounded by hyperplanes perpendicular
to the line $\overline{pq}$ and respectively containing
$p$ and $q$.

\begin{center}
\resizebox{!}{1in}{\includegraphics{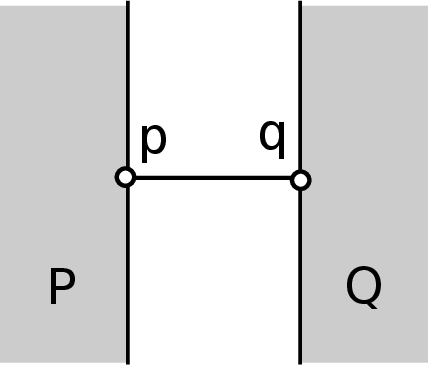}}
\end{center}
\begin{center}
{\bf Figure 3.1:\/} The points $p,q$ and the halfspaces $P,Q$
\end{center}

If $S \cap P=\emptyset$ then
$\R^n-P$ is a convex subset that contains $S$.
This forces $H \subset \R^n-P$, contradicting
$p \in P \cap H$.    Hence
there exists $p' \in S \cap P$.  Likewise
there  exists $q' \in S \cap Q$.  By construction
$\|p'-q'\| \geq \|p-q\|$.
\endproof

Let $F$ and $G$ be the
two boundary components of $\Omega$.
Let $H={\rm Hull\/}(G)$.

\begin{lemma}
  $F$ intersects $H$.
\end{lemma}

\startproof
Assume $F \cap H=\emptyset$.
Choose some $p \in H$.   For $r \geq 1$ let
$D_r$ denote dilation by a factor of $r$ about $p$.   Let
$F_r=D_r(F)$.
Note that
$G \subset H \subset D_r(H)$ and
$F_r \cap D_r(H)=\emptyset$.
Hence $F_r \cap G = \emptyset$.
For large $r$, the loop $F_r$
lies outside
a sphere that strictly contains $G$, and so $F_r$ and $G$ make
a split link.
Since $r \mapsto F_r$ is a smooth isotopy in the complement of $G$, we see that
$F$ and $G$ make a split link.  This is a contradiction.
\endproof

Let $x \in F \cap H$.
The point $x$ is the endpoint of a bend
whose other endpoint $y$ lies in $G \subset H$.
Bends have length at least $1$, so we have
just produced two points $x,y \in H$ whose
distance is at least $1$.  This means that
${\rm diam\/}(H) \geq 1$.  But, by Lemma
\ref{hull}, we have
${\rm diam\/}(G)={\rm diam\/}(H)$.
Hence ${\rm diam\/}(G) \geq 1$.  Since
$G$ is a loop, this means
that $G$ has length at least $2$.
The length of $G$ is $\lambda$, so we see that
$\lambda \geq 2$.
The case $\lambda=2$ is impossible by the same
fattening argument used in the first proof of
Theorem \ref{one}.

\subsection{Proof of Theorem \ref{three}}

As in the first proof of Theorem \ref{three}, let
$\lambda=2+ \epsilon$. Convergence will again be
handled by using ``nearly'' and $\approx$ in the sense
defined in \S \ref{proof31}.

We retain the notation from the preceding subsection.

Let $v$ be the bend whose
endpoints are $x$ and $y$.    Let $P$ be the
plane through $x$ and perpendicular to $v$.
Let $Z$ be the halfspace bounded by $P$
which does not contain $y$.   Since $x \in H$ there
must be some point $z \in G\cap Z$.
Since $\|y-x\| \geq 1$ we have $\|y-z\| \geq 1$.
The loop $G$ has length nearly $2$ and
$G$ contains both $y$ and $z$.  Hence
$G$ is nearly the bigon whose vertices are
$y$ and $z$.    Moreover
$\|y-x\| \approx \|y-z\| \approx 1$.
If follows that $z$ and $x$ are nearly the same point.
There is also a bend $u$ whose endpoints are
$z$ and some $w \in F$.   A similar argument
shows that $F$ must nearly be the bigon with
vertices $w$ and $z$ and
that $\|w-x\| \approx \|w-z\| \approx 1$.

To summarize, we have shown that $F$ and
$G$ are both nearly unit length bigons which nearly
have $x$ as an endpoint.  Moreover, we have
shown that $F$ and $G$ nearly connect the endpoints
of bends $u$ and $v$, respectively.
From here, the proof has the same endgame as
in \S \ref{proof31}.

\section{Proof of Theorem \ref{two}}
\label{exist}

    Using smooth bump functions we make a U-shaped
    curve $U$ which agrees with parallel line
    segments at either end and joins them in the middle.
    The product $U \times [a,b]$ is an isometrically
    embedded rectangle called a {\it pseudofold\/}.
    Now we simply follow the
        instructions for building one of the models depicted in
        Figure 1.1, using pseudofolds
        to smoothly interpolate between the triangular pieces.
        Compare [{\bf HW\/}].
    
    The new object will be an embedded smooth twisted
    cylinder whose boundary is not quite totally
    geodesic. Figure 4.1 shows schematically what
    this looks like.  Now just trim off the rough
    edges. This gives you a smooth embedded twisted paper
    cylinder with slightly larger $\lambda$. With
    this procedure you can make $\lambda$ as close
    to $2$ as you like.  You can check that the boundaries
    make a Hopf link with paper models.
    (We recommend taping string loops along the
    boundary and then pulling them free.)
    
\begin{center}
\resizebox{!}{2in}{\includegraphics{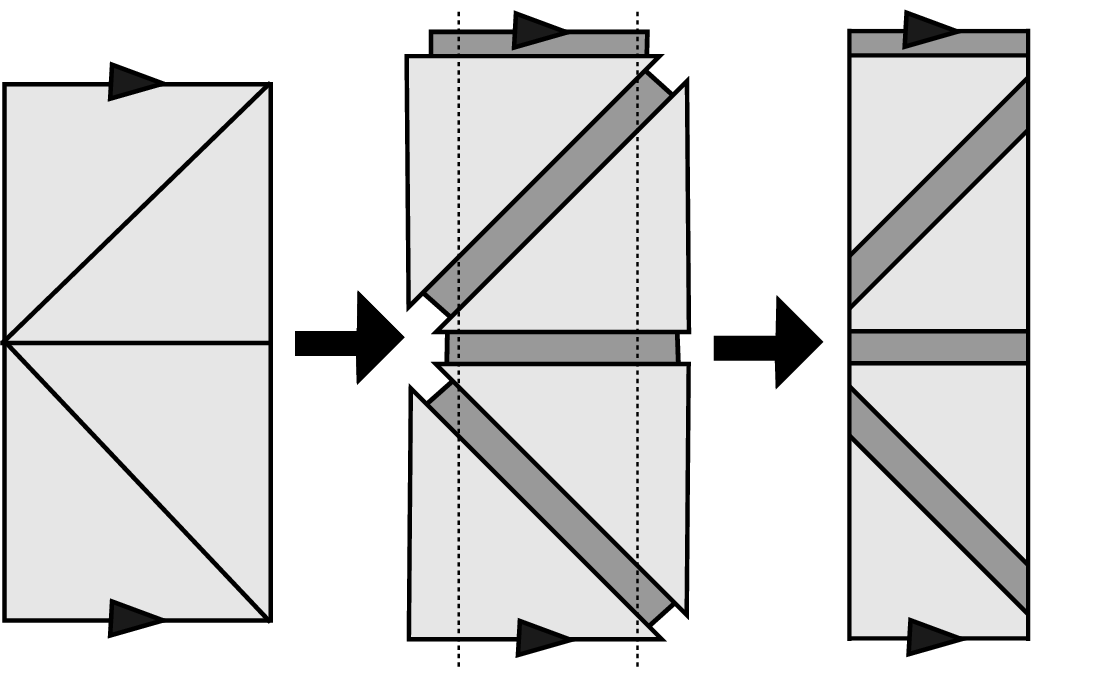}}
\end{center}
\begin{center}
    {\bf Figure 4.1:\/} Perturbing the limiting model
\end{center}
    
    \section{Appendix: The Bend Partition}
\label{app}

Let $I: \Gamma_{\lambda} \to \R^3$ be
an embedded twisted paper cylinder and let
$\Omega=I(\Gamma_{\lambda})$.
  Recall that a bend
  is a straight line segment on $\Omega$ having
  one endpoint in each component of
  $\partial \Omega$. In this appendix we prove
Lemma \ref{BEND}, which  says that
$\Omega$ has a continuous partition into bends.

Lemma \ref{BEND} is closely related to results in classical
differential geometry, which say that developable surfaces
are ruled by straight lines.   Classical treatments tend not to
address certain subtleties covered by Lemma \ref{BEND}, such
as the presence of a boundary and the potential existence of
parts of the surface with zero mean curvature.  This is all
discussed in  [{\bf Spi\/}, Chapter 5, pp. 235--247].

The analogue of Lemma \ref{BEND} is proved e.g. in [{\bf S1\/}] for
the case of Moebius bands and used implicitly
in [{\bf HW\/}].    The proof here is essentially the same.
Let $U \subset \Omega$
denote the nonempty subset where $\Omega$
has nonzero mean curvature.
The following result has many proofs in the literature.
One modern treatment is given in
[{\bf Spi\/}, Corollary 6, p. 241].

\begin{lemma}
  \label{bends}
  $U$ has a partition into bends.
\end{lemma}

Now we derive Lemma \ref{BEND} from Lemma \ref{bends}.
The  bends of $U$ vary continuously because
they are disjoint and all have length at least $1$.
Hence the partition
of $U$ by bends extends to give a continuous
partition
of the closure $\overline U$ by bends.  The
complementary regions of $\Omega-\overline U$
are each flat open trapezoids.  Two opposite parallel
sides of each trapezoid $\tau$ lies in $\partial \Omega$,
and the other two sides are bends of $\overline U-U$.
These bends can be continuously extended to a bend
partition of $\tau$.  Doing this for all
such $\tau$ we extend the bend partition
of $\overline U$ to a bend partition of $\Omega$.
\newline
\newline
{\bf Remark:\/}
In case $I$ is an immersion, Lemma \ref{BEND} would say
that $\Gamma_{\lambda}$ has a partition into prebends.
The proof is the same.  We let $U' \subset \Gamma_{\lambda}$
be the open set of points whose images are points of nonzero
curvature.  Using the same local results we get partition
of $U'$ into prebends, and then we fill in the complementary
trapezoids with prebends.

    \section{References}

\noindent

\noindent
    [{\bf BrS\/}] B. E. Brown and R. E. Schwartz,
    {\it The crisscross and the cup: Two short $3$-twist paper Moebius bands\/},
    preprint 2023, arXiv:2310.10000
    \vskip 9 pt
    \noindent
  [{\bf D\/}] E. Denne, {\it Folded Ribbon Knots in the Plane\/}, 
    The Encyclopedia of Knot Theory (ed. Colin Adams, Erica Flapan, Allison Henrich, Louis H. Kauffman, Lewis D. Ludwig, Sam Nelson)
    Chapter 88, CRC Press (2021)
    \vskip 9 pt
    \noindent
    [{\bf DL\/}] E. Denne, T. Larsen, {\it Linking number and folded ribbon unknots\/},
    Journal of Knot Theory and Its Ramifications, Vol. 32 No. 1 (2023)
    \vskip 9 pt
    \noindent
[{\bf H\/}] A. Hennessey, {\it Constructing many-twist M\"obius bands
  with small aspect ratios\/}.  arXiv:2401:14639, to appear in
Comptes Rendus
\vskip 9 pt
\noindent
[{\bf HW\/}] B. Halpern and C. Weaver,
{\it Inverting a cylinder through isometric immersions and embeddings\/},
Trans. Am. Math. Soc {\bf 230\/}, pp 41--70 (1977)
\vskip 9 pt
\noindent
    [{\bf S1\/}] R. E. Schwartz, {\it The optimal paper Moebius
      Band\/},
    Annals of Mathematics Vol. 201.1 (2025)
\vskip 9 pt
\noindent
    [{\bf S2\/}] R. E. Schwartz, {\it On nearly optimal paper Moebius
      Band\/},
    Advances in Geometry (2025) to appear
\vskip 9 pt
\noindent
[{\bf Spi\/}] M. Spivak, {\it A Comprehensive Introduction to
  Differential Geometry\/}, Vol 3, Third Edition, Publish or Perish (1999)

\end{document}